\documentclass[11pt,twoside,a4paper]{article}

\usepackage{amsmath} \usepackage{amssymb} \usepackage{amsthm}
\usepackage{mathrsfs}

\usepackage[all]{xy} \usepackage{eucal} \SelectTips {cm}{}
\CompileMatrices

\usepackage{a4}
\newtheorem{thm}{Theorem}[section] \newtheorem{lem}[thm]{Lemma}
\newtheorem{prop}[thm]{Proposition} \newtheorem{cor}[thm]{Corollary}

\theoremstyle{definition}

\theoremstyle{remark} \newtheorem{rmk}[thm]{Remark}

\newenvironment{pro}[1][Proof]{{\it{#1:}} }{\hfill $\square$}
\newenvironment{pro*}[1][Proof]{{\it{#1:}} }{}

\numberwithin{equation}{section}

\newcommand{\mengensymb}[1]{\mathbb{#1}}
\newcommand{\N}{\mengensymb{N}} \newcommand{\Z}{\mengensymb{Z}}
\newcommand{\Q}{\mengensymb{Q}} 
\newcommand{\C}{\mengensymb{C}} \newcommand{\F}{\mengensymb{F}}

 \newcommand{\PP}{\mengensymb{P}}

\newcommand{\G}{\mengensymb{G}} 

\newcommand{\BB}{\mathscr{B}} \newcommand{\OO}{\mathcal{O}}

\newcommand{\ov}[1]{\mbox{${\overline{#1}}$}}
\newcommand{\ci}[1]{\mbox{$\overset{{\scriptscriptstyle \circ}}{#1}$}}
\newcommand{\ket}{{\text{\rm k\'et}}} 
\newcommand{\et}{ {\text{\rm\'et}}} 

\newcommand{\DP}{{\rm DP}}

\newcommand{\Hom}{{\rm Hom}} \newcommand{\GL}{{\rm GL}}
\newcommand{\Out}{{\rm Out}} 
\newcommand{\Aut}{{\rm Aut}} 

\newcommand{\inj}{\hookrightarrow} 

\newcommand{\pil}{\pi_1^{\rm log}}
\newcommand{\pita}{\pi_1^{\rm tame}}
\newcommand{\RD}{{\rm R}} \newcommand{\Ho}{{\rm H}}
\newcommand{\Gal}{{\rm Gal}}

\DeclareMathOperator{\Spec}{Spec}

\newcommand{\ph}{\varphi} \newcommand{\eps}{\varepsilon}
\newcommand{\pr}{{\rm pr}} 
\newcommand{\gp}{{\rm gp}} \newcommand{\fs}{{\rm fs}}
\newcommand{\ab}{{\rm ab}}

\newcommand{\sh}{{\rm sh}}

 \newcommand{\red}{{\rm red}}
\newcommand{\lcm}{{\rm lcm}}
\DeclareMathOperator{\rk}{rk} \DeclareMathOperator{\cha}{char}

\begin{document}

\title{\bf A logarithmic view towards \\
  semistable reduction} 
\author{Jakob Stix\thanks{\sc Mathematisches Institut, Universit\"at  Bonn, Beringstra\ss e 1, 53115 Bonn \newline \hspace*{0.45cm} E-mail address: {\tt stix@math.uni-bonn.de}}} 
\date{  } 

\maketitle

\begin{abstract}
  A family of proper, smooth curves yields a monodromy action of the
  fundamental group of the base on the $\Ho^1$ of a fibre. The
  geometric condition of T.~Saito for the action of the wild inertia
  to be trivial is transformed to the condition of logarithmic smooth
  reduction. The proof emphasizes methods and results from logarithmic
  geometry. It applies to quasi--projective smooth curves with \'etale
  boundary divisor.
\end{abstract}


\section{Introduction} \label{secIntro}

The geometry of the degenerate fibres of a proper, generically smooth
curve $f: X \to S$ is very much controlled by the monodromy action on
$\ell$-adic cohomology.  More precisely, let $U \subset S$ be the
smooth locus of $f$ and $\ov{u} \in U$ a geometric point. The sheaf
$\RD^1f_{\ast,\et} \Q_\ell$ is an \'etale local system above $U$,
hence corresponds to an action of $\pi_1(U,\ov{u})$ on the cohomology
$\Ho^1\big(X_{\ov{u}},\Q_\ell\big)$ of the fibre $X_{\ov{u}}$. Let $s$
be a point in the degeneration locus $S - U$ that is a normal point of
codimension 1 in $S$.  Then the monodromy action of the respective
inertia subgroup $I_s \subset \pi_1(U,\ov{u})$, defined up to
conjugation, is nontrivial but quasi-unipotent. It is unipotent if and
only if $X|_U$ admits semistable reduction in $s$.

Nevertheless, trivial action of inertia corresponds to good reduction,
i.e., smoothness. That raises the question about consequences for the
geometry of the reduction if only the wild inertia subgroup acts
trivially.

Let us introduce some notation that remain valid throughout the paper.
We define a {\bf trait} $S$ to be the spectrum of an excellent, strict
henselian discrete valuation ring $R$ with perfect residue field $k$,
uniformizer $\pi$ and field of fractions $K$. We will denote the
generic point by $\eta$ and the closed/special point by $s$. We fix
geometric points $\ov{\eta},\ov{s}$ above $\eta$ and $s$, such that
--- with $\eta^t$ being the maximal tamely ramified subextension of
$\ov{\eta}/\eta$ --- the normalization $S^t$ of $S$ in $\eta^t$ has
$\ov{s}$ as its closed point. These data fix an {\bf inertia subgroup}
$I_s < G_K = \Gal(\ov{K}/K)$ in the absolute Galois group of $K$.  In
case the residue characteristic $\cha(k)=p$ is positive there is also
the {\bf wild inertia subgroup} $P \subset I_s$ which is the $p$-Sylow
subgroup of $I_s$.

A {\bf proper curve} $X/S$ is a flat, proper map $X \to S$ of finite
presentation and pure relative dimension $1$ with geometrically
connected fibres.

On the above question the following answer exists.

\begin{thm}[Saito, {\cite[Theorem 3]{Saito}}] \label{thmSaito}
  Let $S$ be a trait with algebraically closed residue field and
  residue characteristic $p>0$. Let $X/S$ be a proper curve of genus
  $\geq 2$ with smooth generic fibre $X_\eta$ that is minimal with
  respect to $X$ being regular and the reduced special fibre
  $X_{s,\red}$ being a normal crossing divisor on $X$. Then the
  following are equivalent:
\begin{itemize}
\item[(a)] $P$ acts trivially on $\Ho^1_{
    \et}\big(X_{\ov{\eta}},\Q_\ell\big)$ for some prime number $\ell
  \not= p$.
\item[(b)] Any component of the special fibre $X_s$ with multiplicity
  divisible by $p$ is isomorphic to $\PP_k^1$ and intersects the rest
  of the special fibre in exactly two points lying on components of
  multiplicity prime to $p$.
\end{itemize}
\end{thm}

The first impetus for writing this paper was to stress that Saito's
theorem is actually about good reduction --- but in logarithmic
geometry. For an introduction to logarithmic geometry see
\cite{OV,KaJHU}.

The trait $S$ admits a {\bf canonical} fs-log structure: $M_S = \OO_S
- \{0\} \subset \OO_S$. Logarithmic smoothness over $S$ generalizes
semistable maps, i.e., maps that are \'etale locally isomorphic to $\Spec R[t_1,
\dots, t_n]/t_1 \cdot \ldots \cdot t_r - \pi$ over $S$.

In parallel to the classical situation, trivial action of wild inertia
corresponds to an action of the fundamental group of the base
$\pil(S,\ov{\eta}) \cong G_K/P$.  We propose to prove the following
theorem.

\begin{thm} \label{thmlsr}
  Let $S$ be a trait with canonical log structure and residue
  characteristic $p>0$.  Let $X_\eta$ be a proper, log-smooth curve
  over the generic point $\eta \in S$ with negative Euler
  characteristic $\chi = \sum (-1)^q \dim
  \Ho_{\ket}^q\big(X_{\ov{\eta}},\F_\ell\big)$.  Then the following
  are equivalent:
\begin{itemize}
\item[(c)] $P$ acts trivially on
  $\Ho^1_{\ket}\big(X_{\ov{\eta}},\F_\ell\big)$ for some prime number
  $\ell \not= p$.
\item[(d)] $X_\eta$ has good reduction over $S$, i.e., there is a
  proper, log-smooth $X/S$ such that the generic fibre is isomorphic
  to $X_\eta$.
\end{itemize}
Addendum: (1) Moreover, if the above conditions are satisfied, then
the minimal proper, regular model $X/S$ of $X_\eta$ with respect to
the reduced special fibre plus the locus of non-trivial log structure
in the generic fibre being a normal crossing divisor can be endowed
uniquely with a fs-log structure such that $X/S$ is log-smooth.

(2) If there is log-smooth reduction then there is a
log-smooth model (not ne\-ces\-sarily regular) such that no component of
the special fibre has multiplicity divisible by $p$.
\end{thm}

We replaced the coefficients $\Q_\ell$ by $\F_\ell$ because the
statement is apparently stronger. But actually both variants are
equivalent, see Corollary \ref{corcprime}(1).

One should notice, however, that though we might have log-smooth
reduction we still can have non-trivial action from the wild inertia
group on the $\Ho^1\big(X_{\ov{\eta}},\Q_\ell\big)$. For example, look
at an elliptic curve over a $p$-adic field $K$ with semistable
reduction and Tate-model $\C_p^\ast/q^\Z$ such that $q$ is no $p$th
power in $K$.

\subsection*{A guide through the paper}
When the generic fibre has trivial log structure, the fastest way to
Theorem \ref{thmlsr} is via $\S$\ref{secBaseextension}. There one
finds a quick proof relying on the theorem of semistable reduction.
However, $\S\S$\ref{secComputation}--\ref{secCombinatorics} supply a
selfcontained proof based on logarithmic geometry, thus reproving the
semistable reduction theorem. In $\S$\ref{secGroups} and
$\S$\ref{secMultiplicities} the connection between Theorem
\ref{thmSaito} and Theorem \ref{thmlsr} is dealt with.

One purpose of writing these notes was to separate cohomology
($\S$\ref{secComputation}), group theory ($\S$\ref{secGroups}),
logarithmic geometry ($\S$\ref{secSpecialization},
$\S$\ref{secBaseextension}, $\S$\ref{secMultiplicities}) and the
combinatorial argument ($\S$\ref{secCombinatorics}). In particular,
due to the use of logarithmic geometry, the author considers the
combinatorial treatment as easier than in \cite{Saito} or
\cite{Abbes}.

\subsection*{Thanks} I am grateful to Ivan Fesenko for inquiring  about a logarithmic treatment of Saito's theorem. Thanks go to Moshe Jarden for an invitation to Tel Aviv and the excellent working conditions there. Special thanks I owe to Tina \& Ken and of course Sabine for usage of their computers during my stay in KL.

\setcounter{tocdepth}{1} {\scriptsize \tableofcontents }


\section{Computation of cohomology} \label{secComputation}

The basics of logarithmic geometry can be found in \cite{KaJHU},
whereas the fundamentals of Kummer \'etale cohomology are contained in
\cite{OV}. Of this section, only Theorem \ref{thmchi} is used lateron.

We basically encounter fs-log structures of the following type.  Let
$D$ be a divisor on the normal noetherian scheme $X$ with complement
$j:U \to X$.  We denote by $M(\log D)$ the fs-log structure $j_\ast
\G_{m,U} \cap \OO_X \inj \OO_X$ and call it `induced by $D$'.
Sometimes we need to distinguish between the log scheme $X$ and its
underlying scheme $\ci{X}$ and exploit the `forget log' map $\eps : X
\to \ci{X}$.

\subsection{Log-smooth curves} 
We first describe log-smooth curves $X/K$ over a field $K$ with
trivial log structure on $\Spec K$, i.e., $\ci{X}$ is geometrically
connected of pure dimension $1$ and $X \to \Spec K$ is log-smooth.

\begin{lem}
  Let $X/K$ be a log-smooth curve.  Then the underlying scheme
  $\ci{X}/K$ is a (classically) smooth curve and the log structure is
  induced by a divisor $D$ that is relatively \'etale over $K$.
\end{lem}
\begin{pro}
  Let $x \in X$ be a point.  By \cite[3.13]{KaJHU} the map
  $\Omega^1_{X/K} \otimes \kappa(x) \to \ov{M}_{X,x}^{\gp} \otimes
  \kappa(x)$ is surjective.  Hence $\ov{M}_{X,x}^{\gp}$ is of rank
  $\leq 1$.
  
  At points of the open set where the log structure is trivial the
  assertion follows from \cite[Prop 3.8]{KaJHU}.  On the other hand,
  if $\ov{M}_{X,x}^{\gp}$ is of rank $1$, we may choose a generator
  $t$ of the log structure, such that $\Omega^1_{X,x} = \OO_{X,x}
  \cdot d\log t$. By \cite[Prop 3.12]{KaJHU} and \cite[Prop
  3.8]{KaJHU} the corresponding map $X \to \Spec K[\N]$ sending $1 \in
  \N$ to $t$ is \'etale near $x$ in the classical sense. That proves
  the lemma with $D = \{ x \in X \, | \, \rk_\Z \ov{M}_{X,x}^{\gp} = 1
  \}$.
\end{pro}

\subsection{Kummer \'etale cohomology}
Let $X/K$ be a proper, log-smooth curve with $M_X = M(\log D)$.
Kummer \'etale covers of $X$ correspond to covers of $\ci{X}$ that are
at most tamely ramified at most over $D$, hence an isomorphism
$\pil(X) \cong \pita(X,D)$. The usual description of $\Z/\ell$-torsors
shows
\begin{equation} \label{eqTorsors}
\Ho^1_{\ket}\big(X,\Z/\ell\big) =
\Hom\big(\pil(X),\Z/\ell\big) \ .
\end{equation}
Let now $K$ be algebraically closed, and let $\Lambda$ be a finite
ring with $\#\Lambda$ invertible in $K$.  For the computation of the
Kummer \'etale cohomology of $X$ we exploit the spectral sequence
associated to $\eps: X_{\ket} \to \ci{X}_{{\et}}$.  Indeed, via the
logarithmic Kummer sequence we deduce $\big(\Lambda^q\ov{M}_X^\gp\big)
\otimes \Lambda(-q) \cong \RD^q\eps_\ast\Lambda$, see \cite[5.2]{OV},
and thus the following 5-term exact sequence (all coefficients are
constant $\Lambda$):
\begin{equation} \label{eqSequence}
 0 \to \Ho^1_{\et}\big(X\big) \to \Ho^1_{\ket}\big(X\big) \to
 \bigoplus_{x \in D} \ov{M}_{X,x}^\gp \otimes \Lambda(-1)  \xrightarrow{\sum} 
\Ho^2_{\et}\big(X\big) \to \Ho^2_{\ket}\big(X\big) \to 0 \ . 
\end{equation}
We summarize the computation of $\Ho^\ast_\ket\big(X,\Lambda\big)$ by
the following:
\begin{prop} \label{propc}
  All cohomology groups $\Ho_\ket^q(X,\Lambda)$ are free
  $\Lambda$-modules of finite rank, the logarithmic Euler
  characteristic being
  $$\chi(X_{\ket},\Lambda) = \sum_{q} (-1)^q \rk_\Lambda
  \Ho^q_\ket\big(X,\Lambda\big) = 2 - 2g - \deg D $$
  where $g$ is the
  genus of $\ci{X}/K$ and $M_X = M(\log D)$. \hfill $\square$
\end{prop}

\subsection{Logarithmic Euler characteristic} \label{setup}
We will need a combinatorial expression for the Euler characteristic
of a degenerate fibre of a relative curve in terms of the intersection
configuration.  Saito obtained such an expression via inspection of
the vanishing cycles sheave.  We will use Kummer \'etale cohomology
instead.

For intersection theory in this context see for example \cite{Lich}.
In particular, we use intersection numbers $\big(C \bullet D\big)$ on
a regular surface between arbitrary divisors $D$ via their associated
line bundles or more generally line bundles itself and divisors $C$
whose support is proper over some fixed field.

Let $S$ be a trait with its canonical log structure $M(\log s)$. Let
$X/S$ be a proper curve such that the generic fibre $X_\eta$ is
log-smooth, the underlying scheme $\ci{X}$ is regular and the log
structure $M_X$ is induced by a normal crossing divisor $Y_0 + H
\subset X$ where the horizontal part $H$ is finite and generically
\'etale over $S$, and the vertical part $Y_0 = \sum C$ is the divisor
associated to the reduced special fibre $X_{s,\red}$. Here a divisor
$D \subset X$ is called {\bf normal crossing} if \'etale locally $D =
\{t_1 \cdot \ldots \cdot t_r = 0\}$ and $t_1, \dots, t_r \in
\OO_{X,x}$ is part of a regular parameter system, i.e., a regular
system with $t_i$ linear independent in
$\mathfrak{m}_x/\mathfrak{m}_x^2$.

Let $Y_1 = \sum_C m_C C$ be the divisor on $X$ where the sum is taken
over all irreducible components of the special fibre $X_s$ and $m_C=
\#\Z[\frac{1}{p}]/(r_C)$ is the prime-to-$p$ part of the multiplicity
$r_C$ of $C$ in $X_s$ which as a divisor is $Y = \sum r_C C$.

\begin{thm}[compare {\cite[Prop 1.1]{Saito}}] \label{thmchi}
  Let $\Lambda$ be a finite ring such that $\#\Lambda$ is invertible
  in $\OO_S$. The respective logarithmic Euler characteristic of the
  fibres is as follows.
\begin{itemize}
\item[(g)] $\chi\big((X_{\ov{\eta}})_\ket,\Lambda\big) = \sum (-1)^q
  \rk_\Lambda \Ho^q_\ket\big(X_{\ov{\eta}},\Lambda\big) = - \big(Y
  \bullet Y_0 + H + \omega \big)$.
\item[(s)] $\chi\big((X_{\tilde{s}})_\ket,\Lambda\big) = \sum (-1)^q
  \rk_\Lambda \Ho^q_\ket\big(X_{\tilde{s}},\Lambda\big) = - \big(Y_1
  \bullet Y_0 + H + \omega \big)$.
\end{itemize}
Here $\omega$ is the relative dualizing sheaf for $\ci{X}/\ci{S}$, and
$X_{\tilde{s}}$ is the log-geometric fibre in a log-geometric point
$\tilde{s}$ with center in the closed point $s$.
\end{thm}
\begin{pro}
  (g) We apply Grothendieck duality relative $S$ to the coherent sheaf
  $\OO_{X_s}$ on $X$.  As $\omega$ is a line bundle we get the usual
  {\it Adjunction Formula}
  $$
  2\chi\big((X_s)_{\rm ZAR},\OO_{X_s}\big) = - \big(Y \bullet Y +
  \omega \big) \ . $$
  Moreover, by coherent Euler characteristic being
  constant in proper, flat families we have
  $$
  2 - 2g = 2\chi\big((X_\eta)_{\rm ZAR},\OO_{X_\eta}\big)=
  2\chi\big((X_s)_{\rm ZAR},\OO_{X_s}\big) \ , $$
  hence
  $$
  2-2g = - \big(Y \bullet Y + \omega \big) = - \big(Y \bullet Y_0 +
  \omega \big) $$
  as the divisor $Y$ is contained in the kernel of the
  intersection pairing on the special fibre.  Furthermore, $\big(Y
  \bullet H\big) = \deg H/S$ and (g) is reduced to Proposition
  \ref{propc}.
  
  (s) We may assume that the residue field $k$ of the trait $S$ is
  algebraically closed. We have to describe the log-geometric special
  fibre. Let $\tilde{s}$ be the log-geometric point $\tilde{s} =
  \varinjlim_{p \nmid n} s_n$ where $s_n$ is $\Spec k$ with log
  structure $\frac{1}{n}\N \to k$ sending $\frac{1}{n}$ to $0$ and the
  transition maps in the limit are the natural ones. Let $X_n = X_s
  \times_s^\fs s_n$, such that the log-geometric fibre is
  $X_{\tilde{s}} = \varprojlim_{p \nmid n} X_n$.
  
  The set of {\bf double points} of $X_s$ is $\DP(X_s) = \{x \, | \,
  \rk\ov{M}^\gp_{X_s,x} = 2\}$. It consists of the singular locus of
  the normal crossing divisor $Y_0 + H$ and contains $H \cap Y_0$
  ({\bf tails}) and intersections of components of $Y_0$ ({\bf actual
    double points}).  We define for $n \in \N$ the set of double
  points of $X_n$ as $\DP(X_n) = \{x \, | \, \rk \ov{M}^\gp_{X_n,x} =
  2\}$.

\begin{lem}  \label{lemFibregeometry}
  Let $N=\lcm \{m_C\}$ be the prime-to-p part of the least common
  multiple of the multiplicities of components of the special fibre.
\begin{itemize}
\item[(1)] The natural maps $\ci{X}_{\tilde{s}} \to \ci{X}_{nN} \to
  \ci{X}_N$ are bijective closed immersions.
\item[(2)] The natural map $\pr: \ci{X}_{N,\red} \to \ci{X}_{s,\red}$
  is a finite ramified cover at most tamely ramified over the double
  points.  For an irreducible component $C$ of $X_{s,\red}$ we have
  $\sum_{C' \mapsto C} \deg C'/C = m_C$ where the sum is taken over
  components of $\ci{X}_{N,\red}$ mapping to $C$.
\item[(3)] At $x \in C \cap D$ for irreducible components $C,D$ of
  $X_s$ the ramification index for $x' \in C' \mapsto x \in C$ is
  $m_C/\gcd(m_C,m_D)$.
\end{itemize}
\end{lem}
\begin{pro} For (1) see \cite[Section I.3]{Vith}. Up to taking the strict reduced subscheme, the map $s_N \to s$ is a Kummer \'etale cover. The same holds for its fs-base change. Thus the ramification of $ \ci{X}_{N,\red} \to \ci{X}_{s,\red}$ is at most tame over the double points. At a regular point $x$ of a component $C$ the degree (as map of schemes) may be calculated as the cardinality of prime-to-$p$ torsion in 
  $$
  \ov{M}^\gp_{X_s,x} \oplus_{\ov{M}^\gp_{s}} \ov{M}^\gp_{s_N} \cong
  \frac{1}{r_C}\Z \oplus_\Z \frac{1}{N}\Z \cong \Z/m_C \oplus
  \frac{1}{\lcm\{r_C,N\}}\Z \ . $$
  This proves (2). A more careful
  \'etale local analysis at the double points shows (3) which we won't
  use in the sequel.
\end{pro}

To compute the Euler characteristic of
$\Ho_\ket^\ast\big(X_{\tilde{s}},\Lambda\big) = \varinjlim_{p \nmid n}
\Ho_\ket^\ast\big(X_n,\Lambda\big)$ we use the Leray spectral
sequences for the `forget log' maps $\eps_n: X_n \to \ci{X}_n$
\begin{equation} \label{eqSF}
\Ho_\et^p\big(\ci{X}_n,\RD^q\eps_{n,\ast}\Lambda\big) \Longrightarrow
\Ho_\ket^{p+q}\big(X_n,\Lambda\big) 
\end{equation}
which is compatible with the direct limit in $n$. For simplicity we
assume that all $n$ are divisible by $N$ (as in the Lemma
\ref{lemFibregeometry}) and therefore all higher direct images live on
the same \'etale site $\big(\ci{X}_N\big)_\et$.
\begin{lem} The limits of the higher direct images are as follows:
\begin{itemize}
\item[(1)] $ \varinjlim_{p \nmid n} \eps_{nN,\ast} \Lambda = \Lambda$,
\item[(2)] $ \varinjlim_{p \nmid n} \RD^1\eps_{nN,\ast} \Lambda =
  \varinjlim_{p \nmid n} \ov{M}^\gp_{X_{nN}} \otimes \Lambda(-1) \cong
  {\displaystyle \bigoplus_{x \in \DP(X_N)} i_{x,\ast}\Lambda(-1)}$,
\item[(3)] $ \varinjlim_{p \nmid n} \RD^2\eps_{nN,\ast} \Lambda =
  \varinjlim_{p \nmid n} \Lambda^2\ov{M}^\gp_{X_{nN}} \otimes
  \Lambda(-1) = (0)$.
\end{itemize}
\end{lem}
\begin{pro} It suffices to prove (2). Furthermore, by \cite[Thm 5.2]{OV}, we need only to calculate $ \varinjlim_{p \nmid n} \ov{M}^\gp_{X_{nN}}$.
  At points $x$ where $\ov{M}^\gp_{X_N}$ is of rank $1$ we have
  $$
  \varinjlim_{p \nmid n} \ov{M}^\gp_{X_{nN},x} = \varinjlim_{p
    \nmid n} \big(\frac{1}{r_C}\Z \oplus_\Z
  \frac{1}{nN}\Z\big)/\text{\rm tors} = \varinjlim_{p \nmid n}
  \frac{1}{\lcm(r_C,nN)}\Z = \Z_{(p)} . $$
  At a double points $x \in
  \DP(X_N)$ we have
  $$
  \varinjlim_{p \nmid n} \ov{M}^\gp_{X_{nN},x} = \varinjlim_{p
    \nmid n} \left(\big(\frac{1}{r_C}\Z \oplus \frac{1}{r_D}\Z
    \big)\oplus_{\text{\rm diag},\Z} \frac{1}{nN}\Z\right)/\text{\rm
    tors} = \Z \oplus \Z_{(p)} . $$
  Here $\oplus_\Z$ is the coproduct
  over $\Z$ in abelian groups, and $\Z_{(p)}$ is the localization at
  the prime ideal $(p)$.
\end{pro}

Thus the only non-vanishing $E_2$-terms in (\ref{eqSF}) in the limit
over $p \nmid n$ are $E_2^{p0} = \Ho^p_\et\big(\ci{X}_N, \Lambda\big)$
and $E_2^{01} = \bigoplus_{x \in \DP(X_N)} \Lambda(-1)$.

Let $U'= C' - \DP(X_N)$ be the complement of the set of double points
in an irreducible component $C'$ of $X_N$. Let $U$ be the respective
image of $U'$ under $\pr:X_N \to X_s$. The multiplicativity of the
Euler characteristic in tame extensions (Hurwitz formula) the above
yields the following:
\begin{eqnarray*}
 \chi\big((X_{\tilde{s}})_\ket,\Lambda\big) & = &
 \chi\big((\ci{X}_N)_\et,\Lambda\big) - \#\DP(X_n) \\
& = & \sum_{C'} \chi(U'_{\et},\Lambda)
=  \sum_{C} \sum_{C' \mapsto C} \deg(C'/C) \cdot \chi(U_{\et},\Lambda) \\
& = & \sum_{C} - m_C \cdot \big(C \bullet Y_0 + H + \omega\big)
 = - \big(Y_1 \bullet Y_0 + H + \omega \big) \ . 
\end{eqnarray*}
This proves (s).
\end{pro}


\section{$p$-groups acting on $\ell$-groups} \label{secGroups}

This section exploits only elementary group theory.  We denote the
maximal pro-$\ell$ quotient of a pro-finite group $G$ by $G^\ell$.

Let $L$ be a pro-$\ell$ group. Any action on $L$ induces an action on
the maximal abelian $\ell$-elementary quotient $L^\ab/\ell$, which is
an $\F_\ell$ vector space.

\begin{lem} \label{lemPonL}
  A continuous action of a pro-$p$ group on a finitely generated
  pro-$\ell$ group $L$ factors through a finite quotient that maps
  isomorphically onto a subgroup of $\GL(L^{\ab}/\ell)$.
\end{lem}
\begin{pro} 
  This follows from the case of finite groups \cite[Thm 12.2.2]{Hall}.
\end{pro}
\begin{cor} \label{corcprime}
  (1) The monodromy action of the wild inertia in Theorem
  \ref{thmSaito} factors through a finite quotient.
  
  (2) Condition (c) of Theorem \ref{thmlsr} is equivalent to the
  following:
\begin{itemize}
\item[(c')] the restriction to $P$ of the natural $\ell$-adic exterior
  action
  
  $G_K \to \Out\big(\pil(X_{\ov{\eta}})^\ell\big)$ is trivial for some
  $p \not= \ell$.
\end{itemize}
Furthermore, if we are in the situation of Theorem \ref{thmSaito},
condition (a) is equivalent to (c) and (c').
\end{cor}
\begin{pro} (1) The action of $P$ on 
  $\Ho_{\ket}^1\big(X_{\ov{\eta}},\Q_\ell \big)$ comes by scalar
  extension from an action on $\Ho_{\ket}^1\big(X_{\ov{\eta}},\Z_\ell
  \big)$ which is a finitely generated pro-$\ell$ group with
  $\Ho_{\ket}^1\big(X_{\ov{\eta}},\F_\ell \big)$ as maximal
  $\ell$-elementary abelian quotient.  (We use $\Ho^\ast_\ket$ which
  in case of trivial log structure as in Theorem \ref{thmSaito}
  coincides with $\Ho^\ast_\et$ because the argument also applies to
  the situation of Theorem \ref{thmlsr}.)
  
  (2) The group $\Ho_{\ket}^1\big(X_{\ov{\eta}},\F_\ell \big)$ is dual
  to the maximal $\ell$-elementary abelian quotient of
  $\pil(X_{\ov{\eta}})^\ell$. Using a $p$-Sylow subgroup of
  $\Aut\big(\pil(X_{\ov{\eta}})^\ell\big)$ we lift the exterior action
  of $P$ to a true action without changing the image.
\end{pro}
\begin{lem}
  Let $\rho : P \to \Out(L)$ be the exterior action of a pro-$p$ group
  $P$ on a pro-$\ell$ group $L$ corresponding to an extension $1 \to L
  \to G \to P \to 1$.
  
  Then $\rho$ is trivial if and only if $L$ is isomorphic to $G^\ell$
  via the canonical map $L \to G^\ell$.
\end{lem}
\begin{pro}
  Use a $p$-Sylow subgroup of $G$ to split the sequence (from the
  right). $L$ is isomorphic to $G^\ell$ if and only if the sequence
  also admits a retraction $G \to L$ if and only if the action $\rho$
  is trivial.
\end{pro}
\begin{cor} \label{cordescenttotame}
  Condition (c') above is equivalent to the natural map
  $\pil\big(X_{\ov{\eta}}\big)^\ell \to \pil\big(X_{\eta^t}\big)^\ell$
  inducing an isomorphism.
\end{cor}
\begin{pro} The exterior action of $P$ comes from the natural short exact sequence
  $$
  1 \to \pil\big(X_{\ov{\eta}}\big)^\ell \to
  \pil\big(X_{\eta^t}\big)/\ker\Big(\pil\big(X_{\ov{\eta}}\big) \to
  \pil\big(X_{\ov{\eta}}\big)^\ell \Big) \to P \to 1 \ . $$
\end{pro}


\section{Specialization and log-smoothness} \label{secSpecialization}

\subsection{The `easy direction'}
In this section we use results about the logarithmic specialization
map of logarithmic fundamental groups to prove the `easy direction' of
Theorem \ref{thmlsr}, that is `(d) implies (c)'.

\begin{pro} {(Vidal)} (d) implies (c). Let $X/S$ be a log-smooth, proper curve.
  By \cite[Thm I.2.2]{Vith} the logarithmic specialization map induces
  an isomorphism
  $$
  \pil\big(X_{\ov{\eta}}\big)^\ell \to
  \pil\big(X_{\tilde{s}}\big)^\ell$$
  on pro-$\ell$ completions. Hence
  the $G_K$-action factors through $\pil(S) \cong G_K/P$.  Thus
  condition (c') is satisfied, and we are done by Corollary
  \ref{corcprime}.
\end{pro}
\begin{rmk} 
  Alternatively, we may also argue with the vanishing of the sheaf of
  log-vanishing cycles for a log-smooth, proper $X/S$, see \cite[Thm
  3.2]{Nacycles}. Indeed, the sheaf $\RD\Psi^{\log}(\Q_\ell)$ of
  nearby cycles being quasi-isomorphic to $\Q_\ell$ in this case, the
  spectral sequence of nearby cycles
  $$
  E_2^{p,q} =
  \Ho_\ket^p\big(X_{\tilde{s}},\RD^q\Psi^{\log}(\Q_\ell)\big)
  \Longrightarrow \Ho_\ket^{p+q}\big(X_{\ov{\eta}},\Q_\ell\big) $$
  degenerates. Hence the Galois action on
  $\Ho_\ket^{p}\big(X_{\ov{\eta}},\Q_\ell\big) \cong
  \Ho_\ket^p\big(X_{\tilde{s}},\Q_\ell\big)$ factors through $\pil(S)
  \cong G_K/P$.
\end{rmk}

\subsection{Purity}
In fact, the basic ingredient \cite[Thm I.2.2]{Vith} of the proof
above, like its classical counterpart, requires a logarithmic version
of deformation and algebraization together with a purity assertion.
These also yield the following proposition.
\begin{prop} \label{propPurity}
  Let $X/S$ be a proper curve, $H \subset X$ a relative effective
  Cartier divisor. Assume that $X$ is regular and $X_{s,\red} + H$ is
  a normal crossing divisor on $X$. Endow $X$ with the log-regular
  fs-log structure $M(\log X_{s,\red} + H)$. Then the natural maps
  induce isomorphisms
  $$
  \pil\big(X_{\eta^t}\big)^\ell \xrightarrow{\pi_1(j)}
  \pil\big(X_{S^t}\big)^\ell \xleftarrow{\pi_1{(i)}}
  \pil\big(X_{\tilde{s}}\big)^\ell \ . $$
\end{prop}
\begin{pro}  A close inspection of I.~Vidal's proof of \cite[Thm I.2.2]{Vith} shows that log-regular instead of log-smooth over $S$ is sufficient for the deformation and algebraization part. Thus $\pi_1(i)$ is an isomorphism. 
  
  Again, log-regularity is sufficient for the purity result of
  Fujiwara--Kato, cf.\ \cite[Thm 7.6]{OV}, and thus $\pi_1(j)$ is an
  isomorphism.
\end{pro}


\section{Tame base extensions} \label{secBaseextension}

The strategy for the proof of a theorem about good reduction is
composed of three steps: (i) determine a good candidate to work with,
(ii) analyze the good candidate's special fibre and come up with a
cohomological condition that decides whether it is `good', and (iii)
find methods to finally enforce the cohomological condition.

\subsection{The good candidate} \label{secCandidate}
In our case, the good candidate is provided by the theory of minimal
models of surfaces together with a divisor. Let $S$ be a trait with
its canonical log structure $M(\log s)$. As in Section \ref{setup},
let $X/S$ be a proper curve such that the generic fibre $X_\eta$ is
log-smooth, the underlying scheme $\ci{X}$ is regular and the log
structure $M_X$ is induced by a normal crossing divisor $Y_0 + H
\subset X$ where the horizontal part $H$ is finite and generically
\'etale over $S$, and the vertical part $Y_0 = \sum C$ is the divisor
associated to the reduced special fibre $X_{s,\red}$.  Moreover, for
$X/S$ to be our good candidate we ask $X/S$ to be relatively minimal
subject to the above conditions.

\subsection{The fs-base change}
Let $S' \to S$ be a finite tame extension. When equipped with the
canonical log structures it is finite log-\'etale. We denote the
fs-log base change $X \times_S^\fs S'$ by $X'$.  By \cite[Thm
3.3]{Bauer} we have:
\begin{equation} \label{eqbasechange}
 X/S \ \text{ is log-smooth } \ \iff \ X'/S' \ \text{ is log-smooth} . 
\end{equation}
It remains to describe $X'$ to use the above as a `devissage'
argument.  But $X'$ is log-regular itself being log-\'etale over the
log-regular $X$. Thus $\ci{X'}$ is normal \cite[Theorem 4.1]{Kats} and
$M_{X'} = M(\log X'_{s,\red} + H')$ with $H'$ being the preimage of
$H$ under the projection $X' \to X$, see \cite[Theorem 11.2]{Kats}. So
$X'$ is the normalization of $\ci{X} \times_{\ci{S}} \ci{S'}$ and its
singular set is contained in the set where the log structure is of
rank $2$. But $X'$ is not too far away from being the good candidate
over $S'$. We will remedy the deficiencies subsequently.

\subsection{Desingularisation by log blow-ups}
The (classical) desingularisation of a log-regular log scheme can be
performed through subdividing the associated fans like with toric
varieties, see \cite[Section 10]{Kats}. In our case these fans are
two-dimensional. Thus the subdivision is achieved by consecutive log
blow-ups along ideals of $M_X$ generated by two elements. This is of
importance as it implies that there is a log-\'etale desingularisation
map $X'_{{\rm desing}} \to X'$ such that the (reduced) fibres of
positive dimension are chains of projective lines (log blow-ups are
log-\'etale).  Moreover, $X'_{{\rm desing}}$ is again log-regular. On
a log-regular log-scheme whose underlying scheme is regular, the log
structure is induced by a normal crossing divisor.

\subsection{Contraction of $(-1)$-curves}
The underlying scheme of $X'_{{\rm desing}}$ is regular but need not
be relatively minimal. There might exist some $(-1)$-curves that we
want to contract using Castelnuovo's Criterion. But only those
$(-1)$-curves will be contracted that intersect the rest of the
divisor inducing the log structure in at most two points. In fact,
there are no hairs, i.e., $(-1)$-curves that intersect only once, for
otherwise $X/S$ would not have been relatively minimal. On the other
hand, the contraction is also a classical blow-up map whose structure
is under full control. We easily see that it is the underlying map of
a log blow-up. Therefore we obtain a log-\'etale contraction map
$X'_{{\rm desing}} \to X'_{\rm min}$ such that the (reduced) fibres of
positive dimension are chains of projective lines and $X'_{\rm min}$
is relatively minimal with respect to the usual requirements. So there
emerges the good candidate $X'_{\rm min}$ over $S'$. By applying
\cite[Thm 3.3]{Bauer} twice, we have:
\begin{equation} \label{eqdesing}
 X'/S' \ \text{ is log-smooth } \ \iff \ X'_{\rm min}/S' \ \text{ is log-smooth} . 
\end{equation}
Combining (\ref{eqbasechange}) with (\ref{eqdesing}) completes the
devissage argument. While working with our good candidate, we may
perform finite tame extensions of our base trait $S$ without affecting
log-smoothness.

\subsection{A shortcut} \label{subsecShortcut}
Using the theorem of semistable reduction we now get a quick proof of
Theorem \ref{thmlsr} in case $X_\eta$ has trivial log structure ($H =
0$).

\begin{pro} (c) implies (d): The monodromy action of the inertia group on $\Ho_\et^1\big(X_{\ov{\eta}},\Q_\ell\big)$ is known to be quasi-unipotent. As $P$ acts trivial by condition (c) the action is unipotent after a finite tame extension of $S$. The devissage argument above allows to assume that the monodromy action of the inertia group is already unipotent over $S$. 
  
  Then the good candidate $X/S$ is semistable by the theorem of
  semistable reduction in its precise form, cf.\ \cite[Thm 1]{Saito}.
  Hence \'etale locally $X/S$ has the form $\Spec R[t]$ or $\Spec
  R[t_1,t_2]/(t_1\cdot t_2 - \pi)$ and is therefore log-smooth, cf.\ 
  \cite[Ex 3.7(2)]{KaJHU}. So (d) holds for the good candidate.
\end{pro}

The argument above obscures that via logarithmic geometry one may
actually reprove the theorem of semistable reduction. The subsequent
section will achieve this.


\section{Combinatorics within the special fibre} \label{secCombinatorics}

\subsection{A combinatorial argument}
Let $X/S$ be the good candidate of Section \ref{secCandidate}.  Let
$\BB$ be the set of irreducible components of the special fibre.  As
in Section \ref{setup} we introduce the divisors $Y_0 = \sum_{C \in
  \BB} C$, $Y_1 = \sum_{C \in \BB} m_C C$ and $Y = \sum_{C \in \BB}
r_C C$. Remember that $m_C$ was the prime-to-$p$ part of $r_C$, the
multiplicity of $C$ in the special fibre $X_s$.  As above, the
relative dualizing sheaf of $\ci{X}/\ci{S}$ is denoted by $\omega$.
Furthermore, we define the $\Z$-valued linear function $F(D) = \big(D
\bullet Y_0 + H + \omega\big)$ on divisors $D$ with support in the
special fibre.

\begin{thm}[compare {\cite[Section 2]{Saito}}] \label{thmCombinatorics}
  Let $X/S$ be a good candidate as above such that the horizontal part
  $H$ of the divisor inducing the log structure on $X$ is \'etale over
  $S$. If the Euler characteristic $\chi = - \frac{1}{2}F(Y)$ is
  negative, then
  $$F(Y) = F(Y_0) \qquad \text{if and only if} \qquad Y = Y_0 \ .
  $$
\end{thm}
\begin{pro} Clearly we assume that $F(Y) = F(Y_0)$ and have to show that all multiplicities $r_C$ equal $1$.
  
  {\it step 1.} If $Y = r_C \cdot C$ then the theorem follows from
  $F(Y) \not= 0$. We may therefore assume that $\#\BB \geq 2$.
  
  {\it step 2.} Now we locate the non-positive contributions to
  $F(Y)$.
\begin{lem}[{\cite[Lemma 2.4]{Abbes}}] \label{lemTypes}
  Let $\#\BB \geq 2$ and $C \in \BB$. Then $F(C) \leq 0$ occurs if and
  only if $\big(C \bullet C + \omega\big) = - 2$ with $C$ of one of
  the following types:
\begin{itemize}
\item[($\alpha$)] $F(C) = -1$, $C$ intersects only one component $C'$
  of $Y_0$ and avoids $H$. Moreover, $\big(C \bullet C'\big) = 1$ and
  $\big(C\bullet C\big)r_C + r_{C'} = 0$.
\item[($\beta$)] $F(C) = 0$ and $C$ intersects only two components
  $C',C''$ of $Y_0+H$. Moreover, $\big(C \bullet C'\big) = \big(C
  \bullet C''\big) = 1$ and $C'$ has support in $Y_0$. If $C''$ has
  support in $Y_0$ we have $\big(C \bullet C\big)r_C + r_{C'} +r_{C''}
  = 0$ while $\big(C \bullet C\big)r_C + r_{C'} = 0$ if $C''$ has
  support in $H$.
\item[($\gamma$)] $F(C) = 0$, $C$ intersects only one component $C'$
  of $Y_0$ and avoids $H$. Moreover, $\big(C \bullet C'\big) = 2$ and
  $\big(C\bullet C\big)r_C + 2r_{C'} = 0$.
\end{itemize}
\end{lem}
\begin{pro}
  By the {\it Adjunction Formula} $\big(C \bullet C + \omega\big)$
  equals $2g_C-2$ where $g_C$ is the arithmetic genus of $C$. As $g_C
  \geq 0$ and $C$ has to intersect at least one other component of
  $Y_0$ by Zariski connectedness, we have
  $$
  F(C) = \big(C \bullet C + \omega\big) + \sum_{C' \not= C} \big(C
  \bullet C'\big) + \big(C \bullet H\big) \geq -1$$
  Using that $\big(C
  \bullet C + \omega\big)$ is even and $\#\BB \geq 2$ one derives the
  assertion of the lemma. The equation for the multiplicities follows
  from $\big(C \bullet Y\big) = 0$.
\end{pro}

The curves dealt with in Lemma \ref{lemTypes} will be called {\bf of
  type $(\alpha),(\beta)$ or $(\gamma)$} respectively or of {\bf
  non-positive type}.  Due to relative minimality, in all cases
($\alpha$), ($\beta$) or ($\gamma$) the selfintersection is
$\big(C\bullet C\big) \leq -2$.

{\it step 3.} To balance the negative contributions to $F(Y)$ we build
clusters coming from the combinatorial structure of the dual graph of
the special fibre.

First, not all components can be of non-positive type for otherwise
$F(Y)$ is non-positive contradicting the assumption on the Euler
characteristic. Secondly, we define an {\bf $\alpha$-hair} as a
maximal connected subgraph of the dual graph of $Y_0$ that contains
one curve of type $(\alpha)$ and all other curves are of type
$(\beta)$:
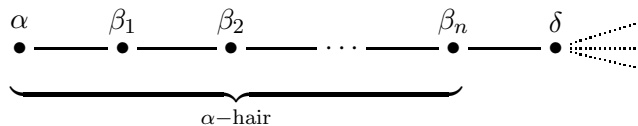
$$
\xymatrix@R=-3pt{{\alpha} & {\beta_1} & {\beta_2} & & {\beta_n} & {\delta} & \\
  {\bullet}
  \ar@{-}[r] & {\bullet} \ar@{-}[r] & {\bullet} \ar@{-}[r] & {\cdots} \ar@{-}[r] & {\bullet} \ar@{-}[r] & {\bullet} \ar@{.}[r] \ar@{.}[ur] \ar@{.}[dr] & \\
  {\phantom{\alpha}} &&&&&&} \vspace{-0.5cm}
$$
$$
\underbrace{ \hspace{6cm} }_{\alpha-\text{hair}} \hspace{2.4cm}
$$
The other neighbour $\delta$ of $\beta_n$ is called the {\bf
  neighbour of the $\alpha$-hair}. The neighour has support in $Y_0$
and is not of non-positive type, for otherwise all components of $Y_0$
were of non-positive type.

Let $\BB_\alpha$ be the set of components that belong to some
$\alpha$-hair, let $\BB_\delta$ be the set of all neighbours of
$\alpha$-hairs, and let $\BB_\rho$ be the set of all other components.
Then $\BB = \BB_\alpha \cup \BB_\delta \cup \BB_\rho$ yields a
disjoint division of the set of irreducible components of $Y_0$.

Let $\pi : \BB_\alpha \to \BB_\delta$ map a component to the neighbour
of the $\alpha$-hair that it is part of. Furthermore, for $D \in
\BB_\delta$ let
$$\ph(D) = (r_D-1)F(D) + \sum_{\{C \in \BB_\alpha \, | \, \pi(C) = D
  \}} (r_C - 1)F(C)$$
be the part in $F(Y) - F(Y_0)$ that is
`combinatorially related' to $D$.
\begin{lem}
  Let $D$ be a neighbour. If $F(Y) > 0$ then $\ph(D) > 0$.
\end{lem}
\begin{pro}
  Let us write $R = Y_0 + H - D - \sum_{\{ C \in \BB_\alpha \, | \,
    \pi(C) = D\}} C$.  Then
  $$
  \ph(D) = \frac{r_D - 1}{2}\Big(F(D) + \big(D \bullet D +
  \omega\big) + \big(D \bullet R\big)\Big) + \! \! \!
  \sum_{\genfrac{}{}{0pt}{}{ C \text{ of type } (\alpha),}{\pi(C) =
      D}} \! \! \!  \frac{r_D - 1}{2} - (r_C - 1) \ . $$
  The following
  lemma shows that the second summand is positive. So if we assume
  that $\ph(D) \leq 0$ then $F(D) + \big(D \bullet D + \omega\big) +
  \big(D \bullet R\big)$ has to be negative. But the minimal values
  are $F(D) = 1$, $\big(D \bullet D + \omega\big) = - 2$, and $\big(D
  \bullet R\big) = 0$ which therefore must be attained.  Hence the
  fibre consists only of $D$ together with at least three
  $\alpha$-hairs neighbouring $D$, for otherwise $D$ were of type
  $(\alpha)$ or $(\beta)$.  Consequently,
  $$F(Y) = \ph(D) + F(Y_0) \leq F(Y_0) = F(D) -
  \#\{\alpha\text{-hairs}\} \leq -2$$
  leads to a contradiction.
\end{pro}
\begin{lem} \label{lemConvex}
  (1) For an $\alpha$-hair together with its neighbour the
  multiplicity $r_C$ is a convex function on the distance from the
  terminal curve of type $(\alpha)$.
  
  (2) Let $r_\alpha$, resp. $r_\delta$, be the multiplicity of the
  terminal curve of type $(\alpha)$, resp. the neighbour, of the same
  $\alpha$-hair. The convex function of (1) is strict monotone
  increasing and in particular $r_\alpha \leq \frac{1}{2}r_\delta$.
\end{lem}
\begin{pro}
  This follows from the relations between the multiplicities and the
  bound on the selfintersection numbers listed in Lemma
  \ref{lemTypes}.
\end{pro}

{\it step 4.} We claim that under the assumptions of $F(Y) =F(Y_0)$
there is no curve of type $(\alpha)$. Indeed, we have only
non-negative summands in
$$
0 = F(Y)-F(Y_0) = \sum_{D \in \BB_\delta} \phi(D) + \sum_{C \in
  \BB_\rho} (r_C - 1)F(C) $$
forcing all summands to vanish. Hence all
components which are neither of type $(\beta)$ nor $(\gamma)$ are
reduced and there do not exist any $\alpha$-hairs by step 3.

{\it step 5}. It remains to deal with curves of type $(\beta)$ or
$(\gamma)$.  As not all components are of non-positive type the curves
of type $(\beta)$ organize in maximal chains in the dual graph of
$Y_0$ and no chain forms a closed circle. At both ends these chains
terminate: either in a curve which is not of non-positive type and
thus is reduced, or in a component $C$ of type ($\beta$) which
intersects $H$. Because $H/S$ is \'etale, a local equation for $C$ at
$C \cap H$ is given by $\pi \cdot \text{\rm unit}$ and thus $r_C = 1$.
An adaptation of Lemma \ref{lemConvex} shows that the multiplicity,
being a convex function on the position of the component in the chain,
must be a constant equal to $1$ along such a $\beta$-chain.

A curve $C$ of type $(\gamma)$ intersects with a curve $C'$ which is
not of non-po\-sitive type.  Thus $\big(C \bullet C\big)r_C + 2r_{C'}
= 0$, and $r_{C'} = 1$ shows $r_C = 1$ for not to contradict relative
minimality.  This completes the proof of Theorem
\ref{thmCombinatorics}.
\end{pro}

\subsection{Proof of the main theorem}
We may now complete the proof of Theorem \ref{thmlsr} apart from the
special properties of the special fibres in the case of logarithmic
smooth reduction. Those will be discussed in $\S$\ref{secMultiplicities}.

\begin{pro} (c) implies (d). Let $X/S$ be the good candidate as in Section \ref{secCandidate}. 
  By Corollary \ref{corcprime}, Corollary \ref{cordescenttotame} and
  Proposition \ref{propPurity}, condition (c) implies that
  $\pil(X_{\ov{\eta}})^\ell \cong \pil(X_{\tilde{s}})^\ell$ and via
  (\ref{eqTorsors}) thus $\chi\big((X_{\ov{\eta}})_\ket,\Z/\ell\big) =
  \chi\big((X_{\tilde{s}})_\ket,\Z/\ell\big)$ for that particular
  prime $\ell \not= p$.
  
  From the 5-term exact sequence (\ref{eqSequence}) we obtain a Galois
  sequence
  $$
  \Ho^1_\ket\big(X_{\ov{\eta}},\Z_\ell\big) \to
  \Z_\ell(-1)[H(\ov{K})] \to \Z_\ell(-1) . $$
  Trivial action of
  $P$ on $\Ho^1_\ket$ yields that $H(\ov{K})=H(K^{\rm sep})$ is
  $P$-stable.  Thus $H/S$ is at most tamely ramified and $H$ splits
  after an apropriate finite tame extension of the base trait $S$.

  The d\'evissage argument of Section \ref{secBaseextension} allows us
  to assume that $H/S$ is \'etale and that the support of $Y_1-Y_0$ is
  contained in chains of projective lines, i.e., $\PP^1$'s that
  intersect the rest of the fibre in two points. In particular $F(Y_1)
  = F(Y_0)$.  By Theorem \ref{thmchi} we conclude further that
  $$
  F(Y) = -\chi\big((X_{\ov{\eta}})_\ket,\Z/\ell\big) =
  -\chi\big((X_{\tilde{s}})_\ket,\Z/\ell\big) = F(Y_1) = F(Y_0) \ . $$
  The combinatorial reasoning of Theorem \ref{thmCombinatorics} shows
  that $Y=Y_0$ or that the special fibre is a reduced normal crossing
  divisor that meets the relative \'etale $H$ with normal crossing.
  It remains to argue that such an $X/S$ is actually log-smooth.
  
  At regular points of $Y_0+H$ the map $X \to S$ is strict and even
  smooth in the classical sense ($k$ is perfect). At a double point
  $x$, we choose regular parameters $u,v$ such that $\{u \cdot v = 0\}
  = Y_0 + H$ locally at $x$. We may impose $uv^{a} = \pi$ with either
  $a=1$ (actual double point) or $a=0$ (tail). In both cases $X/S$
  factors \'etale locally as
  $$
  \Spec \OO_{X,x}^\sh \xrightarrow{j} W = \Spec R[u,v]/(uv^a = \pi)
  \xrightarrow{h} S$$
  with log structure on $W$ induced by $u^\N
  v^\N$. The map $h$ is log-smooth being the fs-base change of a map
  induced by $(1,a): \N \to \N^2$. The map $j$ is essentially \'etale
  as can be checked on completions which are regular rings of
  dimension $2$. Indeed, by the choice of $u,v$ the map $j$ is
  surjective on the Zariski cotangent space after scalar extension.
\end{pro} 
\begin{cor}[semistable reduction of curves]
  Let $S$ be a trait. Any proper, smooth curve of genus $\geq 2$ over
  the generic point $\eta \in S$ admits a semistable model over some
  finite extension of $S$.
\end{cor}
\begin{pro}
  By Corollary \ref{corcprime} the wild inertia action vanishes after
  some finite extension. Thus we are reduced to the above proof of (c)
  implies (d).
\end{pro}


\section{Multiplicities divisible by $p$} \label{secMultiplicities}

In this section we address the addendum to Theorem \ref{thmlsr} that
describes the special fibre in the case of log-smooth reduction. We
emphasize geo\-me\-tric reasoning thereby providing a link between the
conditions (b) and (d).

Let $X/S$ be a proper, log-smooth curve such that $\ci{X}$ is
relatively minimal with respect to being regular and the divisor
inducing the log structure being normal crossing. 

Let $E$ be a component of $X_s$ with multiplicity divisible by $p$.
The reasoning of Sections \ref{secBaseextension} \&
\ref{secCombinatorics} shows that $E$ is a projective line that
intersects the rest of the special fibre in exactly two points. It
remains to show that if $x$ is in the intersection of two components
$C,D$ of the special fibre, then not both components have multiplicity
divisible by $p$. This follows from the fact that by \cite[Thm
3.5]{KaJHU} the cokernel of the map $\ov{M}^\gp_{S,s} \to
\ov{M}^\gp_{X,x}$ has no $p$-torsion. Indeed, this cokernel is
isomorphic to $\Z \times \Z/\gcd(r_C,r_D)$. Thus (d) implies (b).

For the converse that (b) implies (d) we contract an annoying
projective line $E$ with multiplicity divisible by $p$ and
intersecting $C,D$ transversally. The resulting normal $S$-curve
$X'/S$ has a rational singularity in $x'$, the image of $E$. With the
help of M.~Artin's work on rational singularities we determine the
\'etale local structure: from the monoid $P = \{(c,d,e) \in \N^3 \, |
\, c+d+\big(E \bullet E\big)e \geq 0 \}$ with $(r_C,r_D,r_E) = p \in
P$ we form the log-smooth $S$-scheme $\Spec R[P]/\pi=p$ that is
\'etale locally isomorphic at $P=0$ to $X'$ in $x'$, see
\cite[I.3.4.3]{Stix}.  After the contraction, the curve in general is
not regular any more. Never\-the\-less, it is log-smooth over $S$ thus
proving the addendum of Theorem \ref{thmlsr}.


\end{document}